\newcommand{\dfrac}{\displaystyle\frac}

\documentclass[12pt]{article}

\usepackage{times}
\usepackage{bm}
\usepackage{natbib}
\usepackage{graphicx}

\begin{document}

\title{On relationship between regression models and interpretation of multiple regression coefficients}

\author{A. N. Varaksin, V. G. Panov}

\maketitle

\begin{abstract}
In this paper, we consider the problem of treating linear regression equation coefficients in the case of correlated predictors. It is shown that in general there are no natural ways of interpreting these coefficients similar to the case of single predictor. Nevertheless we suggest linear transformations of predictors, reducing multiple regression to a simple one and retaining the coefficient at variable of interest. The new variable can be treated as the part of the old variable that has no linear statistical dependence on other presented variables.
\end{abstract}


\section{Introduction }

 Regression analysis is one of the main methods for studying dependency factors in diverse fields of inquiry where use of statistical  methods is expedient (see e.g. \citet{DS98}). The efficiency of its application depends on the model and the set of explanatory variables (predictors) chosen. The most popular regression model is described by a linear equation expressing the dependence of the mean value  of the variable (response, outcome) to be explained on the set of predictors.

 The natural applicability domain of regression analysis is a case of continuous outcome and predictors. In this area, the classical  regression analysis theory provides a thorough description of outcome dependence on explanatory variables considered. Of most  interest in the linear model are the coefficients at predictors. For example, in the simplest case of single predictor $X_{1} $ and  dependent variable $Y$ the linear regression equation is given by

\[
y=b_{0} +b_{1} x_{1}
\]

Factor $b_{1} $ is proportional to the coefficient of correlation between response $Y$ and predictor $X_{1} $. Furthermore, $b_{1} $ represents an increase (or a decrease, if $b_{1} $ is negative) in the mean of $Y$ associated with a 1-unit increase in the value of $X,\ X = x + 1$ versus $X = x$. The sign of $b_{1} $ indicates the trend in the relationship between $Y$ and $X_1$.

 Such thorough information about the relationship between outcome and single explanatory variable makes one wish to treat the  coefficients of a multiple regression equation in a similar manner. It is well known, however, that in linear multiple regression  models such interpretation of regression coefficients is not correct if there are correlations among predictors (\citet{DS98}, \citet{Nalimov}, \citet{Ehrenberg}). Moreover, in some practical cases such interpretation is in conflict with common sense (\citet{Varaksin}, see below section~\ref{app2realData}). The unique case where interpretation of multiple regression equation coefficients is meaningful is pairwise  statistical independence of predictors. Then multiple regression coefficients coincide with corresponding simple regression  coefficients for the outcome on a particular predictor (\citet{DS98}).

 Thus, the presence of correlated predictors renders the identification of the biomedical meaning of multiple regression equation  coefficients a difficult task. Association among predictors or among predictors and outcome leads to unpredictable changes in regression coefficients and results in a loss of meaning in each particular coefficient.

 Nevertheless we cannot confine ourselves to independent (uncorrelated) variables only, as in most applications of regression analysis  there are important problems with correlated predictors, e.g. various air pollution rates (see bellow section \ref{app2realData}). Another important example is  epidemiological studies (research into disease prevalence and its association with risk factors). Such factors as sex and age are  invariably present in epidemiological data, being related to both other independent variables and outcome. These inherent variables  which confuse the effect on the response and other predictors are called confounders. Taking into account confounders in
 data analysis presents a difficult problem that does not have any correct solution as yet.

 In a range of biomedical applications of regression analysis, of major interest is some variable $X_{1} $ which is considered along  with accompanying variables $X_{2} ,X_{3} ,...,X_{k} $ (confounders). Upon finding a multiple regression equation that depends on all  of these predictors one has to treat coefficient $b_{1} $ standing at the principal predictor, with all other predictors adjusting  the action of main variable $X_{1} $. We shall consider below a way to interpret $b_{1} $ in terms of simple regression of outcome on  a new variable, $X_{1}^{*} $. For simplicity, we shall discuss cases of two and three predictors. The general case may be considered
 in a similar way.

\section{Regression equation with two predictors}

 Let us consider continuous variables $Y,X_{1} .X_{2} $ and corresponding linear regression equation for outcome $Y$ on predictors  $X_{1} .X_{2} $

\begin{equation} \label{f1}
y=b_{0} +b_{1} x_{1} +b_{2} x_{2}
\end{equation}

As usual, we suppose that coefficients $b_{0} ,b_{1} ,b_{2} $ and other regression coefficients below have been obtained by the least squares method. We assume that the (linear) dependence of response $Y$ on predictor $X_{1} $ is significant, so $b_{1} \ne 0$. Finally, let the linear regression equation with response $X_{1} $ and predictor $X_{2} $ be given by

\[x_{1} =c_{120} +c_{12} x_{2} \]

We define a new variable, $X_{1}^{*} $, in which the linear dependence of $X_{1} $ on $X_{2} $ `is excluded' as follows
\[
X_{1}^{*} =X_{1} -c_{12} X_{2}
\]

Let us build a simple regression equation describing the mean of outcome $Y$ as a function of new predictor $X_{1}^{*} $

\begin{equation} \label{f2}
y=a_{10}^{*} +a_{1}^{*} x_{1}^{*}
\end{equation}

We have pairs of corresponding variables: $X_{1} $ and $x_{1} $, $X_{1}^{*} $ and $x_{1}^{*} $. Obviously, these variables cannot be interchanged; in particular, variables $X_{1} $, $X_{1}^{*} $ cannot be substituted in equations (\ref{f1}) and (\ref{f2}) instead of $x_{1} $ and $x_{1}^{*} $, respectively. If it were possible, one might transcribe equations (\ref{f1}) and (\ref{f2}) as

\[
\begin{array}{l}
{y=b_{0} +b_{1} \left(x_{1} +\dfrac{b_{2} }{b_{1} } x_{2} \right)} \\ {y=a_{10}^{*} +a_{1}^{*} \left(x_{1} -c_{12}
x_{2} \right)}
\end{array}
\]

Although these equations are different, they have the same slope, as follows from the following theorem.

{\bf Theorem~1.}\label{theorem1}
\textit{In equation (\ref{f1}), coefficient $b_{1} $ is equal to coefficient $a_{1}^{*} $ in equation (\ref{f2}), i.e.
\begin{equation} \label{f3}
b_{1} =a_{1}^{*} ,
\end{equation}
and it is possible that} $\dfrac{b_{2} }{b_{1} } \ne -c_{12} $.

A similar statement holds for coefficients $b_{2} $ and $a_{2}^{*} $, where $a_{2}^{*} $ is the coefficient at variable $x_{2}^{*} $ in a simple regression equation $y=a_{20}^{*} +a_{2}^{*} x_{2}^{*} $, and a new variable $X_{2}^{*} $ is defined from the regression equation $x_{2} =c_{210} +c_{21} x_{1} $ as $X_{2}^{*} =X_{2} -c_{21} X_{1} $.

Formal proof of Theorem~\ref{theorem1}  is provided in Appendix~1.

Now coefficient $b_{1} $ of multiple regression equation (\ref{f1}) may be treated as follows. Recall that $b_{1} $ cannot be interpreted per se. But it is equal to coefficient $a_{1}^{*} $ of simple regression model (\ref{f2}). Hence we transform the problem of interpretation of $b_{1} $ into one of interpretation of a new variable $X_{1}^{*} $. It is easy to check that $X_{1}^{*} $ and $X_{2} $ are uncorrelated. So one can say that variable $X_{1}^{*} $ is obtained from variable $X_{1} $ by excluding the part of it that is linearly dependent on it. This does not mean that by constructing variable $X_{1}^{*} $ we can split the contributions of $X_{1} $ and $X_{2} $ to response $Y.$ In fact, there is no way to do this given correlated predictors.

 Now consider a more general way to define variable $X_{1}^{*} $, namely, let $X_{1}^{*} =X_{1} -\gamma X_{2} $, where $\gamma  $ is a real number, and pose the question: how many values may $\gamma $ take for equality (\ref{f3}) to hold?  In the case under consideration, we can express the dependence of \textit{$a_{1}^{*} $ }on parameter $\gamma $ in explicit  form as follows

\begin{equation} \label{f4}
a_{1}^{*} \left(\gamma \right)=\dfrac{\overline{X_{1} Y}-\overline{X_{1} }\,\overline{Y}-\gamma \left(\overline{X_{2} Y}-\overline{X_{2}
}\,\overline{Y}\right)}{var\left(X_{1} \right)-2\gamma cov\left(X_{1} ,X_{2} \right)+\gamma ^{2} var\left(X_{2} \right)} ,
\end{equation}
where the bar over a symbol denotes the mean of the variable, \textit{var} and \textit{cov }stand for variance and covariance, respectively.

{\bf Theorem~2.}\label{theorem2}
\textit{Equation $a_{1}^{*} \left(\gamma \right)=b_{1} $ has two solutions, videlicet
\[
\gamma _{1} =c_{12} ,\quad \gamma _{2} =-\dfrac{b_{2} }{b_{1} }
\]
}

 For the proof of this theorem we refer the reader to Appendix 2.

\begin{figure}
\begin{center}
\includegraphics[width=8cm,height=5cm,keepaspectratio]{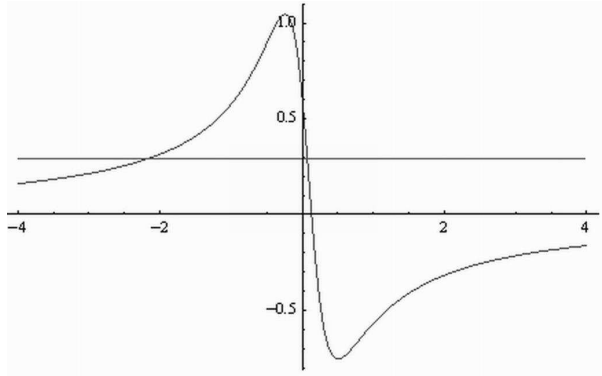}
\caption{Plot of regression coefficient $a_{1}^{*} $ as a function of $\gamma $.}
\label{2Dplot}
\end{center}
\end{figure}

 Given the explicit expression for $a_{1}^{*} \left(\gamma \right)$ in formula (\ref{f4}), we can plot it (see Fig.~\ref{2Dplot},  where artificial data is used with $b_{1} =0.2918$ which is drawn as a horizontal line). There are some general properties in  $a_{1}^{*} \left(\gamma \right)$: it is defined throughout the real axis, has two extrema, and the real axis is an asymptote to it.

\section{Regression equation with three predictors}

 Now consider the case of one outcome $Y$ and three predictors $X_{1} ,X_{2} ,X_{3} $. The point of interest is predictor $X_{1} $ the  other predictors being confounders. We want to get an interpretation of coefficient $b_{1} $ at the variable of interest in the  multiple regression equation

\begin{equation} \label{f5}
y=b_{0} +b_{1} x_{1} +b_{2} x_{2} +b_{3} x_{3}
\end{equation}

 We can introduce the regression equation of $X_{1} $ on covariates $X_{2} ,X_{3} $:

\[
x_{1} =c_{0123} +c_{12} x_{2} +c_{13} x_{3} ,
\]
and define a new variable by the formula

\begin{equation} \label{f6}
X_{1}^{*} =X_{1} -c_{12} X_{2} -c_{13} X_{3}
\end{equation}
As in section 1, we could find a simple regression equation for $Y$on covariate $X_{1}^{*} $

\begin{equation} \label{f7}
y=a_{01}^{*} +a_{1}^{*} x_{1}^{*}
\end{equation}
Similar to Theorem~\ref{theorem1}, we have the following statement.

{\bf Theorem~3.}\label{theorem3}
\textit{
Coefficient $b_{1} $ of equation (\ref{f5})  is equal to coefficient $a_{1}^{*} $ of
equation (\ref{f7}), that is
}
\[
b_{1} =a_{1}^{*}
\]

The proof of this theorem is given in Appendix 3.

\begin{figure}
\begin{center}
\includegraphics[width=300pt,height=250pt,keepaspectratio]{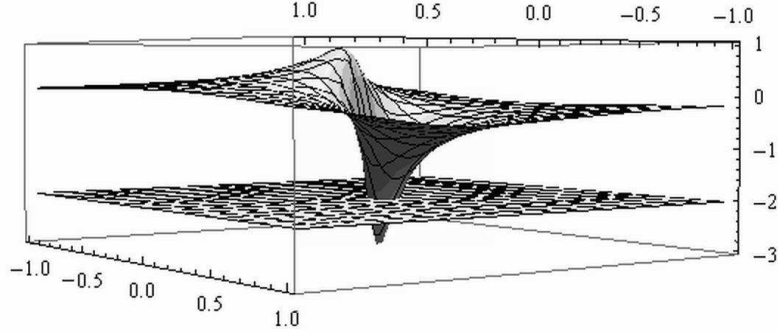}
\caption{Surface $z=a_{1}^{*} \left(\gamma _{2} ,\gamma _{3} \right)$ and plane $z=b_{1} $.}
\label{3Dplot}
\end{center}
\end{figure}

Going over to a more general case, we can define covariate $X_{1}^{*} $ as follows

\[
X_{1}^{*} =X_{1} -\gamma _{2} X_{2} -\gamma _{3} X_{3} ,
\]
where $\gamma _{2} ,\gamma _{3} $ are some real numbers. Then regression coefficient $a_{1}^{*} $ becomes a function of two real variables $\gamma _{2} ,\gamma _{3} $. The shape of surface $z=a_{1}^{*} \left(\gamma _{2} ,\gamma _{3} \right)$ is shown in Figure~\ref{3Dplot} (using simulated data with $b_{1} =-2.031$).

As one can see from Fig.~\ref{2Dplot} and Fig.~\ref{3Dplot}, the character of the dependence of $a_{1}^{*} $ on corresponding parameter(s) in both cases is similar. The same is true of the general case.

\section{Applications to real data analysis}
\label{app2realData}
\subsection{Regression with two predictors}

Let us consider the use of Theorem~\ref{theorem1}  for investigating the dependency of incidence on various air pollution toxicants of City St.-Petersburg (Russia). The primary data were published in \citet{Scherbo}. In the remainder of this section, we assume incidence to be incidence rate in the adult population (i.e. the number of disease cases per 1000 adult population a year) averaged over a 5-year observation period. In the primary data, the rates of incidence were gathered across 19 boroughs of St.-Petersburg. We consider toxicant concentrations as random variables, i.e. mean toxicant concentration expressed in maximum concentration limit (MCL) terms and averaged over 5-year observation period. Each of these variables takes on 19 values in accordance with the number of boroughs. We denote these covariates by the usual chemical notations: $CO, NO_2, SO_2, Pb$ etc. (the data consists of 12 pollutants).

 The simple linear regression equations of response $Y$(incidence) on concentrations of $CO$ and $NO_2$ are given by

\begin{equation} \label{f8}
Y=603+579\, CO
\end{equation}

\begin{equation} \label{f9}
Y=414+416\, NO_{2}
\end{equation}

According to equation (\ref{f8}), incidence increases by 579 cases per 1000 population at an increase in $CO$ concentration by MCL unit a year. Equation (\ref{f9}) may be interpreted in the same way. In short, both $CO$ and $NO_2$ increase incidence.

There is a tight positive correlation between predictors $CO$ and $NO_2$. Pearson's correlation coefficient is 0.75, and the regression equation is

\[
CO=-0.131+0.576\, NO_{2}
\]

This shows that growth in one toxicant is related to growth in another. Hence, one can conjecture that equation (\ref{f8}) does describe an increase in incidence at a simultaneous increase in both pollutants ($CO$ and $NO_2$). A question then arises: could one specify the `pure' influence of each toxicant on incidence, separating the contribution of one toxicant from that of the other?

 To extract the contribution of each toxicant to the incidence in the presence of other toxicants, researchers often use a multiple  regression equation including all toxicants. Such interpretation is common in some biological and medical applications of regression  analysis. We refer to \citet{McNamee} as a typical exposition. In the case under consideration, we obtain a multiple regression equation

\begin{equation} \label{f10}
Y=465+390\, CO+191\, NO_{2}
\end{equation}

  A lot of authors consider the coefficients of a multiple regression equation obtained by means of the least squares method to be meaningless if there are correlations among predictors (\citet{DS98}, \citet{Aivazian}, \citet{Ehrenberg}). These coefficients cannot be used to   assess separately the dependence of $Y$ on $CO$ and $Y$ on $NO_2$. Nevertheless, there are other authors who treat each   coefficient of a multiple regression equation as the contribution of an individual toxicant to the outcome against the background of other toxicants (e.g. \citet{McNamee}). Moreover, this contribution has to be refined as compared to   (\ref{f8})--(\ref{f9}). Their supposition is that predictors as if distribute their influence on the outcome   in a multiple regression equation so that each predictor describes its influence with the other being in the background. According  to this viewpoint, the addition of another toxicant, $NO_2$, to $CO$ and change from (\ref{f8}) to (\ref{f10}) should attenuate the effect of $CO$ because the corresponding coefficient diminished from 579 to 390. The same conclusion holds for $NO_2$ and $CO$ and equations (\ref{f9}) and (\ref{f10}).

 These authors do not provide any substantive explanation for the biomedical meaning of variations in the coefficients in  (\ref{f8})--(\ref{f9}) and (\ref{f10}); nor do they explain the refined contribution of each  individual toxicant. Variations in regression coefficients could be explained by going over from simple regressions  (\ref{f8}) or (\ref{f9}) to multiple regression (\ref{f10}). Indeed, coefficient $b_1=390$ in  equation (\ref{f10}) is equal to coefficient $a_{1}^{*} $ in the simple regression equation

\[Y=a_{10}^{*} +a_{1}^{*} CO^{*} ,\]
where covariate $CO^{*}$ is defined by

\begin{equation} \label{f11}
CO^{*} =CO-0.576\, NO_{2}
\end{equation}

By (\ref{f11}), predictor $CO^{*}$ is obtained from $CO$ by excluding its part correlated with $NO_2$. Then $b_1=390$ means an increased incidence rate at a growth in $CO$ concentration excluding the linear statistical dependence of $ CO$ and $NO_2$.

One can similarly treat coefficient $b_2 = 191$ in (\ref{f10}). It is equal to $a_{2}^{*} $ in the simple regression equation

\[
Y=a_{20}^{*} +a_{2}^{*} \, NO_{2}^{*} ,
\]
where $NO_{2}^{*} $ is a part of toxicant $NO_2$ which contains no linear statistical dependence on $CO.$

We seem to have obtained a consistent picture: by excluding the (linear) dependence of one toxicant on the other we arrive at a `pure' influence of a particular factor on incidence. Since both factors increase the incidence, and the concentration of each factor increases with growth in the other, one can anticipate that the magnitudes of the coefficients in equation (\ref{f10}) should be less than in (\ref{f8})--(\ref{f9}). This is exactly so in the case under consideration.

It is not as simple as that though. Let us consider the dependence of incidence $Y$ on the concentrations of $CO$ and $SO_2$. A simple regression equation of $Y$ on $SO_2$ is given by

\[
Y=919+52\, SO_{2}
\]
The association between $CO$ and $SO_2$ is very similar to that between $CO$ and $NO_2$. For instance, the correlation coefficient is 0.73 and the regression equation is

\begin{equation} \label{f12}
CO=0.272+0.316\, SO_{2}
\end{equation}
The multiple regression equation in the case considered is (\citet{Varaksin})

\begin{equation} \label{f13}
Y=634+1047\, CO-278\, SO_{2}
\end{equation}

Assuming the coefficients of (\ref{f13}) to be refined ones we should treat the magnitude 1047 as a `pure' influence of $CO$ against the background of $SO_2$, and  $-278$ as a `pure' influence of $SO_2$ against the background of $CO.$ Obviously, such interpretation of regression coefficients is invalid, since the `pure' influence of toxicant $SO_2$ becomes negative. The reason for such misinterpretation is the tight correlation between predictors $CO$ and $SO_2$. One has to take into account this correlation in treating regression coefficients.

 The coefficient at $CO$ in (\ref{f13}) is twice as large as that in (\ref{f8}). By Theorem~\ref{theorem1}, coefficient $b_1 = 1047$ is equal to the slope in

\begin{equation} \label{f14}
Y=697+1047\, CO^{*} ,
\end{equation}
where $CO^{*} =CO-0.316\, SO_{2} $. In biomedical terms, we obtain an inexplicable picture: we have reduced the toxic burden on the population by removing one of the two toxicants, but the incidence grows with $CO$ even more rapidly. In mathematical terms, we can explain this as follows. It is clear from the definition of $CO^{*} $ that its range is less than the range of $CO.$ In both cases, the incidence is the same, which implies an increase in coefficient $b_1$.
 Generally, inequality $b_{1} >a_{1} $ is impossible if we consider the multiple regression coefficients as refined ones. But if we  refer to equality (\ref{appB_2}), we can see that under $a_{2} \ll a_{1} $ and correlation coefficient $r$ close to 1, inequality $b_{1} >a_{1} $  may hold true. The formula (\ref{appB_3}) also explains the possibility of a negative value for coefficient $b_{2} $.

\subsection{Regression with three predictors}

 Let us consider a regression equation of incidence $Y$ on three predictors $CO, NO_2$ and $SO_2.$ By the least square  method, we obtain the equation

\[
Y={\rm \; }494{\rm \; }+{\rm \; }857\, CO+{\rm \; }194\, NO_{2} - {\rm \; }279\, SO_{2}
\]
Equation (\ref{f5}) becomes

\[
CO{\rm \; }=-0.108{\rm \; }+{\rm \; }0.386\, NO_{2} +{\rm \; }0.195\, SO_{2} \; \;
\]
The new variable $CO^{*} $ is defined by (\ref{f6}), and the simple regression equation for $Y$ on this predictor is given
by

\[
Y={\rm \; }1076{\rm \; }+{\rm \; }857\, CO^{*}
\]
We see that $b_{1} =a_{1}^{*} $ as well.

 Note that the correlation coefficient of model (\ref{f12}) is $r=0.74$, and that of model (\ref{f14}) is
 $r=0.46$. The latter is less than the coefficient of correlation between incidence $Y$ and $CO$ ($r=0.58$).

\section{Conclusion}

 Let there be two regression equations for outcome $Y$

\[
y=a_{0} +a_{1} x_{1}
\]
and
\[
y=b_{0} +b_{1} x_{1} +b_{2} x_{2}
\]
If predictors $X_{1} ,X_{2} $ are uncorrelated, then $a_{1} =b_{1} $. Hence, inequality $a_{1} \ne b_{1} $ is caused by the presence  of correlation between the predictors. What is the epidemiological meaning of changing coefficient $a_{1} $ to $b_{1} $ after adding predictor $X_{2} $ to the simple regression model? Is the influence of the predictors on the outcome redistributed between them?  The answer is definitely `no'. Usually, the addition of a second covariate is aimed at taking into account the combined effect of predictors on outcome. But what does `take into account' mean? There are no reasonable explanations of this term.

 In view of Theorem~\ref{theorem1}  we can state that the addition of $X_2$ to regression equation $y=b_{0} +b_{1} x_{1} $ brings us to regression  equation $y=a_{10}^{*} +a_{1}^{*} x_{1}^{*} $. The new variable $X_{1}^{*} $ contains no linear statistical dependence on $X_{2} $. A similar interpretation holds for the case of three variables as well as for the general one.

\section*{Acknowledgement}
The research was supported by Basic Science Research for Medicine Program through the Ural Division of RAS Presidium (12-$\Pi$-2-1033).

\appendix
\section*{Appendix 1}
\label{App_A}
\subsection*{Proof of the Theorem~\ref{theorem1} }

Let us first prove a technical statement, being of significance in its own right.
Let there be a set of predictors $X_{1} ,X_{2},...,X_{k} $ and let $Y_0$ be an outcome. The values of $p$ observations over predictors and the outcome combine into matrices $X$ and $Y$

$$
X=\left(\begin{array}{cccc} {1} & {X_{11} } & {...} & {X_{k1} } \\ {1} & {X_{12} } & {...} & {X_{k2} } \\ {\vdots } & {\vdots } &
{\vdots } & {\vdots } \\ {1} & {X_{1p} } & {...} & {X_{kp} } \end{array}\right),\quad Y=\left(\begin{array}{c} {Y_{1} } \\ {Y_{2} } \\
{\vdots } \\ {Y_{p} } \end{array}\right)
$$

The first column contains unities so that we have the same formulae for calculating $b_0$ in the same way as other $b_i.$
Let $B$ denote the column of coefficients $b_0, b_1, b_2 ,\dots , b_k$. To find a linear regression equation for response $Y$ from predictors $X_{1} ,X_{2} ,\dots ,X_{k} $, we have to minimize the mean square residual of $Y$ and $X\,B$ i.e.

\begin{equation}
\mathop{\min }\limits_{\mathbf{B}} (Y-X\,B)(Y-X\,B)^{T}\label{appA_1}
\end{equation}
where the ${}^T$ denotes matrix transposition. The problem (\ref{appA_1}) has a unique solution under the usual least squares method assumptions, e.g. if the matrix $X^TX$ is invertible (see e.g. \citet{DS98}[Chapter 5]). Such assumption will be needed throughout Appendix~1.

 Let $\Gamma $ denote a nonsingular square matrix of order $k$

\[
\Gamma =\left(\begin{array}{cccc} {\gamma _{11} } & {\gamma _{12} } & {...} & {\gamma _{1k} } \\ {...} & {...} & {...} & {...} \\
{\gamma _{k1} } & {\gamma _{k2} } & {...} & {\gamma _{kk} } \end{array}\right)
\]
and let $C$ be a matrix of order $\left(k+1\right)\times \left(k+1\right)$

\[
C=\left(\begin{array}{ccccc} {1} & {0} & {0} & {...} & {0} \\ {0} & {\gamma _{11} } & {\gamma _{12} } & {...} & {\gamma _{1k} } \\
{...} & {...} & {...} & {...} & {...} \\ {0} & {\gamma _{k1} } & {\gamma _{k2} } & {...} & {\gamma _{kk} }
\end{array}\right)=\left(\begin{array}{cc} {1} & {0} \\ {0} & {\Gamma } \end{array}\right)
\]

Let us introduce a vector, $X_{0} =\left(1,X_{1} ,X_{2} ,...,X_{k} \right)$, and consider linear transformation of variables $X_{1}
,X_{2} ,...,X_{k} $ by the matrix $C$

\begin{equation}
X_{0}^{*} =X_{0} \, C\label{appA_2}
\end{equation}

Thus, the new variables $X_{0}^{*} =\left(1,X_{1}^{*} ,X_{2}^{*} ,...,X_{k}^{*} \right)$ obtained from variables $X_{0} =\left(1,X_{1},X_{2} ,...,X_{k} \right)$ by means of linear transformation are given by

\[
X_{i}^{*} =\sum _{j=1}^{k}\gamma _{ij} X_{j}
\]

Finally, we denote by $X^{*} $ a matrix constructed from $X_{0}^{*} $ in the same way as $X$  from $X_0$, and $B^{*} $ stands for the column of coefficients $b_{0}^{*} ,b_{1}^{*} ,...,b_{k}^{*} $.

\textbf{Proposition~1.}
\label{prop_A1}
\textit{Let multiple regression equation of outcome $Y_0$ on predictors $X_{1} ,X_{2} ,...,X_{k} $ be
\[
y=\sum _{i=0}^{k}b_{i} x_{i}
\]
Then coefficients $b_{0}^{*} ,b_{1}^{*} ,...,b_{k}^{*} $ of the multiple regression equation for $Y_0$ on predictors $X_{1}^{*} ,X_{2}^{*} ,...,X_{k}^{*}$
\[
y=\sum _{k=0}^{n}b_{k}^{*} x_{k}^{*}
\]
can be found from the matrix equality}

\[
B^{*} =C^{-1}\, B
\]

It is easy to check that under condition (\ref{appA_2}) we have

\begin{equation} \label{appA_3}
X^{*} =X\,C
\end{equation}

To find a regression equation relative to new variables $X_{i}^{*} $ we need to solve the minimization problem

\[
\mathop{\min }\limits_{B^{*} } (Y-X^{*}\, B^{*} )(Y-X^{*} B^{*} )^T
\]

Given equality (\ref{appA_3}), we have

\[
(Y-X^{*} B^{*} )(Y-X^{*} B^{*} )^T=(Y-X\,C\,B^{*} )(Y-X\,C\,B^{*} )^T
\]

Hence we obtain $C\, B^{*} =B$, since the minimization problem (\ref{appA_1}) has a unique solution. It is obvious from its definition that matrix $C$ is invertible and

\[
C^{-1} =\left(\begin{array}{cc} {1} & {0} \\ {0} & {\Gamma ^{-1} } \end{array}\right)
\]

This brings us to the end of the proof of the Proposition.

\textbf{Proof}
To prove Theorem~\ref{theorem1}, consider the case of two predictors $X_1$, $X_2$, and matrix $C$ is equal to

\[
C=\left(\begin{array}{ccc} {1} & {0} & {0} \\ {0} & {1} & {0} \\ {0} & {-c_{12} } & {1} \end{array}\right)
\]

Then $X_{1}^{*} =X_{1} -c_{12} X_{2} $ and $X_{2}^{*} =X_{2} $. Applying the Proposition to matrix $C,$ we obtain

\begin{equation} \label{appA_4}
B^{*} =\left(\begin{array}{c} {b_{0}^{*} } \\ {b_{1}^{*} } \\ {b_{2}^{*} } \end{array}\right)=C^{-1} B=\left(\begin{array}{ccc} {1} &
{0} & {0} \\ {0} & {1} & {0} \\ {0} & {c_{12} } & {1} \end{array}\right)\left(\begin{array}{c} {b_{0} } \\ {b_{1} } \\ {b_{2} }
\end{array}\right)=\left(\begin{array}{c} {b_{0} } \\ {b_{1} } \\ {b_{2} +c_{12} b_{1} } \end{array}\right)
\end{equation}

 It can be easily seen that $X_{1}^{*} $ and $X_{2} $ are uncorrelated (correlation coefficient is equal to zero). Therefore,
 coefficients $b_{1}^{*} ,b_{2}^{*} $ of the multiple regression equation for outcome $Y$ on predictors $X_{1}^{*} ,X_{2} $

\[
y=b_{0}^{*} +b_{1}^{*} x_{1}^{*} +b_{2}^{*} x_{2}
\]
are equal to the corresponding coefficients of the simple regression equations for $Y$ on predictors  $X_{1}^{*} $ and $X_{2} $, respectively

\begin{equation} \label{appA_5}
\begin{array}{l} {y=a_{01}^{*} +a_{1}^{*} x_{1}^{*} ,\quad y=a_{02} +a_{2} x_{2} } \\ {b_{1}^{*} =a_{1}^{*} ,\quad b_{2}^{*} =a_{2} }
\end{array}
\end{equation}

Using (\ref{appA_4}), we obtain $b_{1}^{*} =b_{1} $, and combining this with (\ref{appA_5}) we obtain $b_{1} =a_{1}^{*}
$, which finishes the proof.

\section*{Appendix 2}
\subsection*{Proof of the Theorem~\ref{theorem2}}

 Let us consider linear transformation of variable $X_{1} $

\begin{equation}
X_{1}^{*} =X_{1} -\gamma X_{2}\label{appB_1}
\end{equation}
Then in the regression equation

\begin{equation} \label{B_f2}
y=a_{10}^{*} +a_{1}^{*} x_{1}^{*}
\end{equation}
coefficient $a_{1}^{*} $ becomes a function of parameter $\gamma $. Its
explicit expression is given by

\begin{equation} \label{B_f4}
a_{1}^{*} \left(\gamma \right)=\dfrac{\overline{X_{1} Y}-\overline{X_{1} }\,\overline{Y}-\gamma \left(\overline{X_{2} Y}-\overline{X_{2}
}\,\overline{Y}\right)}{var\left(X_{1} \right)-2\gamma cov\left(X_{1} ,X_{2} \right)+\gamma ^{2} var\left(X_{2} \right)} ,
\end{equation}

\textbf{Proof} (of Theorem~2)
Recall the following regression equations for outcome $Y$ on predictors $X_{1} ,X_{2} $ (jointly and separately)

\[
\begin{array}{l} {y=b_{0} +b_{1} x_{1} +b_{2} x_{2} } \\ {y=a_{01} +a_{1} x_{1} } \\ {y=a_{02} +a_{2} x_{2} } \end{array}
\]
and we introduce matrices

\[
A=\left(a_{1} ,a_{2} \right),\quad B=\left(b_{1} ,b_{2} \right),\quad C=\left(\begin{array}{l} {\begin{array}{cc} {1} & {c_{12} }
\end{array}} \\ {\begin{array}{cc} {c_{21} } & {1} \end{array}} \end{array}\right),
\]
where $c_{ij}$  are regression coefficients from the equations

\[
\begin{array}{l} {x_{1} =c_{012} +c_{12} x_{2} } \\ {x_{2} =c_{021} +c_{21} x_{1} } \end{array}
\]

According to the theorem \citet{PV10}, we have equality $A=B\cdot C$. Now, suppose that $C$ is an invertible matrix (the
opposite case is discussed below in the Remark\ref{rem1}). Then

\[
B=A\cdot C^{-1}
\]

Thus we obtain the following representation of regression coefficients $b_{1} ,b_{2} $ ($r$ denotes the correlation coefficient
between $X_{1} ,X_{2} $)

\begin{equation}
b_{1} =\dfrac{a_{1} -a_{2} c_{21} }{1-c_{12} c_{21} } =\dfrac{a_{1} -a_{2} c_{21} }{1-r^{2} }\label{appB_2}
\end{equation}

\begin{equation}
b_{2} =\dfrac{a_{2} -a_{1} c_{12} }{1-c_{12} c_{21} } =\dfrac{a_{2} -a_{1} c_{12} }{1-r^{2} }\label{appB_3}
\end{equation}

From (\ref{appB_2})--(\ref{appB_3}), it follows

\begin{equation}
\begin{array}{l} {a_{1} -a_{2} c_{21} =b_{1} \left(1-r^{2} \right),\quad a_{2} -a_{1} c_{12} =b_{2} \left(1-r^{2} \right)} \\ {b_{1}
c_{12} =a_{2} -b_{2} ,\quad b_{2} c_{21} =a_{1} -b_{1} } \\ {c_{12} =\dfrac{a_{2} }{b_{1} } -\dfrac{b_{2} }{b_{1} } =\dfrac{a_{2} -b_{2}
}{b_{1} } ,\quad c_{21} =\dfrac{a_{1} }{b_{2} } -\dfrac{b_{1} }{b_{2} } =\dfrac{a_{1} -b_{1} }{b_{2} } ,\quad } \\ {r^{2} =c_{12} c_{21}
=\dfrac{a_{1} a_{2} }{b_{1} b_{2} } -\dfrac{a_{1} }{b_{1} } -\dfrac{a_{2} }{b_{2} } +1,\quad 1-r^{2} =1-\dfrac{a_{1} }{b_{1} }+\dfrac{a_{2} }{b_{2} } -\dfrac{a_{1} a_{2} }{b_{1} b_{2} } } \end{array}\label{appB_4}
\end{equation}

If we equate the right hand side of (\ref{appB_2}) to the right hand side of (\ref{B_f4}), we obtain the roots of the
equation  $a_{1}^{*} \left(\gamma \right)=b_{1} $ (after some simplification)

\[\begin{array}{l} {\gamma _{1,2} =\dfrac{1}{2(a_{1} -a_{2} c_{21} )var\left(X_{2} \right)} } \\ {\left(\begin{array}{l}
{2cov\left(X_{1} ,X_{2} \right)\left(a_{1} -a_{2} c_{21} \right)-a_{2} var\left(X_{2} \right)\left(1-c_{12} c_{21} \right)\pm } \\
{\sqrt{\begin{array}{l} {4\left(-a_{2} +a_{1} c_{12} \right)c_{21} var\left(X_{1} \right)var\left(X_{2} \right)\left(-a_{1} +a_{2}
c_{21} \right)+} \\ {\left(-2a_{1} cov\left(X_{1} ,X_{2} \right)+a_{2} \left(2c_{21} cov\left(X_{1} ,X_{2} \right)+var\left(X_{2}
\right)\left(1-c_{12} c_{21} \right)\right)\right)^{2} } \end{array}} } \end{array}\right)} \end{array}\]

Applying (\ref{appB_4}), we get

\[\begin{array}{l} {\gamma _{1,2} =\dfrac{1}{2var\left(X_{2} \right)b_{1} \left(1-r^{2} \right)} } \\ {\left(\begin{array}{l} {2r\sigma
\left(X_{1} \right)\sigma \left(X_{2} \right)b_{1} \left(1-r^{2} \right)-a_{2} var\left(X_{2} \right)\left(1-r^{2} \right)\pm } \\
{\sqrt{\begin{array}{l}{4b_{1} b_{2} c_{21} \left(1-r^{2} \right)^{2} var\left(X_{1} \right)var\left(X_{2} \right)+}\\{\left(-2a_{1} r\sigma \left(X_{1}
\right)\sigma \left(X_{2} \right)+a_{2} \left(2c_{21} r\sigma \left(X_{1} \right)\sigma \left(X_{2} \right)+var\left(X_{2}
\right)\left(1-r^{2} \right)\right)\right)^{2} }\end{array}} } \end{array}\right)} \end{array}\]
where $\sigma \left(X_{i} \right)=\sqrt{var\left(X_{i} \right)} $.

Next, we expand the second summand in the radicand and factor out the $-2r\sigma \left(X_{1} \right)\sigma \left(X_{2} \right)$. After that, $a_{1} -a_{2} c_{21} $ is substituted by $b_{1} \left(1-r^{2} \right)$ (see (\ref{appB_4})). We get

\[\begin{array}{l} {\gamma _{1,2} =\dfrac{1}{2var\left(X_{2} \right)b_{1} \left(1-r^{2} \right)} } \\ {\left[\begin{array}{l}
{\left(1-r^{2} \right)\left(2r\sigma \left(X_{1} \right)\sigma \left(X_{2} \right)b_{1} -a_{2} var\left(X_{2} \right)\right)\pm } \\
{\sqrt{4b_{1} b_{2} c_{21} \left(1-r^{2} \right)^{2} var\left(X_{1} \right)var\left(X_{2} \right)+\left(-2r\sigma \left(X_{1}
\right)\sigma \left(X_{2} \right)b_{1} \left(1-r^{2} \right)+a_{2} var\left(X_{2} \right)\left(1-r^{2} \right)\right)^{2} } }
\end{array}\right]} \end{array},\]

or

\[\gamma _{1,2} =\dfrac{1}{2b_{1} } \left[2rb_{1} \dfrac{\sigma \left(X_{1} \right)}{\sigma \left(X_{2} \right)} -a_{2} \pm \sqrt{4b_{1}
b_{2} c_{21} \dfrac{var\left(X_{1} \right)}{var\left(X_{2} \right)} +\left(2rb_{1} \dfrac{\sigma \left(X_{1} \right)}{\sigma \left(X_{2}
\right)} -a_{2} \right)^{2} } \right]\]

Applying (\ref{appB_4}) again, we obtain the required equalities

\[\begin{array}{l} {\gamma _{1,2} =\dfrac{1}{2b_{1} } \left[2b_{1} c_{12} -a_{2} \pm \sqrt{4b_{1} b_{2} c_{12} +\left(2b_{1} c_{12}
-a_{2} \right)^{2} } \right]=}\\ {c_{12} -\dfrac{a_{2} }{2b_{1} } \pm \sqrt{\left(\dfrac{a_{2} }{b_{1} } -c_{12} \right)c_{12} +\left(c_{12}
-\dfrac{a_{2} }{2b_{1} } \right)^{2} } =} \\ {c_{12} -\dfrac{a_{2} }{2b_{1} } \pm \sqrt{\dfrac{a_{2} }{b_{1} } c_{12} -c_{12}^{2}
+c_{12}^{2} -\dfrac{a_{2} }{b_{1} } c_{12} +\left(\dfrac{a_{2} }{2b_{1} } \right)^{2} } =c_{12} -\dfrac{a_{2} }{2b_{1} } \pm
\left|\dfrac{a_{2} }{2b_{1} } \right|} \end{array}\]
That is
\[\gamma _{1} =c_{12} ,\quad \gamma _{2} =c_{12} -\dfrac{a_{2} }{b_{1} } \]
or
\[\gamma _{1} =c_{12} ,\quad \gamma _{2} =-\dfrac{b_{2} }{b_{1} } \]

\textbf{Remark}\label{rem1}
 If correlation matrix $C$ is singular, then $r^{2} =1$, i.e. predictors $X_{1} ,X_{2} $ are proportional. In this case, the problem of finding a multiple regression equation on variables $X_{1} ,X_{2} $ cannot be posed, since it leads to an inconsistent system of linear equations.

\section*{Appendix 3}
\label{App_C}
\subsection*{Proof of the Theorem~\ref{theorem3}}
 The method of proving Theorem~\ref{theorem3}  as considered in this Appendix contains the main ideas of the proof of the general statement.

 Let the linear multiple regression equation for outcome $Y$ on predictors $X_{1} ,X_{2} ,X_{3} $ be

$$ 
y=b_{0} +b_{1} x_{1} +b_{2} x_{2} +b_{3} x_{3}
$$

We introduce new variable $X_{1}^{*} $ by

\[X_{1}^{*} =X_{1} -\gamma _{2} X_{2} -\gamma _{3} X_{3} ,\]
where $\gamma _{2} ,\gamma _{3} $ are some constants. So we perform linear transformation of predictors by the matrix

\[C=\left(\begin{array}{cccc} {1} & {0} & {0} & {0} \\ {0} & {1} & {0} & {0} \\ {0} & {-\gamma _{2} } & {1} & {0} \\ {0} & {-\gamma
_{3} } & {0} & {1} \end{array}\right),\quad \det \left(C\right)=1\]

The inverse matrix $C^{-1} $ is equal to

\[
C^{-1} =\left(\begin{array}{cccc} {1} & {0} & {0} & {0} \\ {0} & {1} & {0} & {0} \\ {0} & {\gamma _{2} } & {1} & {0} \\ {0} &
{\gamma _{3} } & {0} & {1} \end{array}\right)
\]
Hence the coefficients of the regression equation for $Y$ on $X_{1} ,X_{2} ,X_{3} $ and that on $X_{1}^{*} ,X_{2}^{*} =X_{2} ,X_{3}^{*} =X_{3} $ are connected by (see Appendix~1)

\begin{equation} \label{appC_2}
\left(\begin{array}{c} {b_{0}^{*} } \\ {b_{1}^{*} } \\ {b_{2}^{*} } \\ {b_{3}^{*} } \end{array}\right)=\left(\begin{array}{cccc} {1} &
{0} & {0} & {0} \\ {0} & {1} & {0} & {0} \\ {0} & {\gamma _{2} } & {1} & {0} \\ {0} & {\gamma _{3} } & {0} & {1}
\end{array}\right)\cdot \left(\begin{array}{c} {b_{0} } \\ {b_{1} } \\ {b_{2} } \\ {b_{3} } \end{array}\right)=\left(\begin{array}{c}
{b_{0} } \\ {b_{1} } \\ {\gamma _{2} b_{1} +b_{2} } \\ {\gamma _{3} b_{1} +b_{3} } \end{array}\right)
\end{equation}
In particular, for arbitrary $\gamma _{2} ,\gamma _{3} $ coefficients $b_{1} $ and $b_{1}^{*} $ are equal. From now on we assume $\gamma _{2} =c_{12} ,\gamma _{3} =c_{13} $.

Let us consider a simple regression equation for $Y$ on $X_{1}^{*} $

\[y=a_{01}^{*} +a_{1}^{*} x_{1}^{*} \]

We have divided the proof of Theorem~\ref{theorem3} into a sequence of lemmas.

\textbf{Lemma~1.}\label{lemma1}
The multiple correlation coefficient of variable $X_{1}^{*} $ on predictors $X_{2} ,X_{3} $ is equal to zero.

\textbf{Proof}
Let $\rho _{1\cdot 23}^{*} $ be the multiple correlation coefficient of $X_{1}^{*} $ on variables $X_{2} ,X_{3} $. By its definition

\[
\left(\rho _{1\cdot (23)}^{*} \right)^{2} =1-\dfrac{\left|Corr\right|}{C_{11} } ,
\]
where $\left|Corr\right|$ is the determinant of the correlation matrix of variables \textit{$X_{1}^{*} ,X_{2} ,X_{3} $,} $C_{11}$ is the cofactor of the (1,1) entry of the matrix $Corr.$ Therefore

\[
\left|Corr\right|=\left|\begin{array}{ccc} {1} & {r_{12} } & {r_{13} } \\ {r_{21} } & {1} & {r_{23} } \\ {r_{31} } & {r_{32} } & {1}
\end{array}\right|,C_{11} =\left|\begin{array}{cc} {1} & {r_{23} } \\ {r_{32} } & {1} \end{array}\right|
\]
Similar to the case of two predictors, one can see that $r_{12} = r_{21} = r_{13} = r_{31} = 0.$ Hence

\[
\left|Corr\right|=\left|\begin{array}{ccc} {1} & {0} & {0} \\ {0} & {1} & {r_{23} } \\ {0} & {r_{32} } & {1}
\end{array}\right|=C_{11} =\left|\begin{array}{cc} {1} & {r_{23} } \\ {r_{32} } & {1} \end{array}\right|,
\]
i.e. $\rho_{1\cdot{23}}^{*} =0.$

\textbf{Lemma~2.}\label{lemma2}
Let us have multiple regression equations for outcome Y on variables $X_{1}^{*} ,X_{2} $

\begin{equation} \label{appC_3}
y=b'_{0} +b'_{1} x_{1}^{*} +b'_{2} x_{2}
\end{equation}
and that on variables $X_{1}^{*} ,X_{3} $

\[
y=b''_{0} +b''_{1} x_{1}^{*} +b''_{3} x_{3}
\]
Also, consider simple regression equations for outcome Y on predictors $X_{2} $ and $X_{3} $ respectively

\[
\begin{array}{l} {y=a_{02} +a_{2} x_{2} } \\ {y=a_{03} +a_{3} x_{3} } \end{array}
\]
Then

\[
b'_{2} =a_{2} ,b''_{3} =a_{3}
\]
Besides,

\begin{equation} \label{appC_4}
b'_{1} =b''_{1} =a_{1}^{*} ,
\end{equation}
where $a_{1}^{*} $ is the regression coefficient from equation (\ref{f2}).

Proof. As it is mentioned above, covariates $X_{1}^{*} ,X_{2} $ are uncorrelated as well as $X_{1}^{*} ,X_{3} $. Hence $b'_{2} =a_{2},\ b''_{3} =a_{3} $. The last equality (\ref{appC_4}) is implied by Theorem~\ref{theorem1}.

\textbf{Lemma~3.}\label{lemma3}
Let there be a multiple regression equation for variable $X_{3} $ $X_{1}^{*} ,X_{2} $

\begin{equation} \label{appC_5}
x_{3} =\alpha _{0312} +\alpha _{31}^{*} x_{1}^{*} +\alpha _{32} x_{2}
\end{equation}
Then
\[
\alpha _{31}^{*} =0
\]
For the multiple regression equation of $X_{2} $ on predictors $X_{1}^{*} ,X_{3} $
\[
x_{2} =\alpha _{0213} +\alpha _{21}^{*} x_{1}^{*} +\alpha _{23} x_{3}
\]
we have
\[
\alpha _{21}^{*} =0
\]

\textbf{Proof}
We obtain it by a tedious calculation. By the least squares method, coefficient $\alpha _{31}^{*} $ can be obtained from a system of linear equations. The numerator of the expression for $\alpha _{31}^{*} $ is the determinant

\begin{equation} \label{appC_6}
\left|\begin{array}{ccc} {1} & {\overline{X_{3} }} & {\overline{X_{2} }} \\ {\overline{X_{1} }-c_{12} \overline{X_{2} }-c_{13}
\overline{X_{3} }}\quad & {\overline{X_{1} X_{3} }-c_{12} \overline{X_{2} X_{3} }-c_{13} \overline{X_{3}^{2} }}\quad & {\overline{X_{1} X_{2}
}-c_{12} \overline{X_{2}^{2} }-c_{13} \overline{X_{2} X_{3} }} \\ {\overline{X_{2} }} & {\overline{X_{2} X_{3} }} &
{\overline{X_{2}^{2} }} \end{array}\right|
\end{equation}
From corresponding systems of linear equations we obtain

\[c_{12} =\dfrac{\left|\begin{array}{ccc} {1} & {\overline{X_{1} }} & {\overline{X_{3} }} \\ {\overline{X_{2} }} & {\overline{X_{1}
X_{2} }} & {\overline{X_{2} X_{3} }} \\ {\overline{X_{3} }} & {\overline{X_{1} X_{3} }} & {\overline{X_{3}^{2} }}
\end{array}\right|}{\left|\begin{array}{ccc} {1} & {\overline{X_{2} }} & {\overline{X_{3} }} \\ {\overline{X_{2} }} &
{\overline{X_{2}^{2} }} & {\overline{X_{2} X_{3} }} \\ {\overline{X_{3} }} & {\overline{X_{2} X_{3} }} & {\overline{X_{3}^{2} }}
\end{array}\right|} ,\quad c_{13} =\dfrac{\left|\begin{array}{ccc} {1} & {\overline{X_{2} }} & {\overline{X_{1} }} \\ {\overline{X_{2}
}} & {\overline{X_{2}^{2} }} & {\overline{X_{1} X_{2} }} \\ {\overline{X_{3} }} & {\overline{X_{2} X_{3} }} & {\overline{X_{1} X_{3}
}} \end{array}\right|}{\left|\begin{array}{ccc} {1} & {\overline{X_{2} }} & {\overline{X_{3} }} \\ {\overline{X_{2} }} &
{\overline{X_{2}^{2} }} & {\overline{X_{2} X_{3} }} \\ {\overline{X_{3} }} & {\overline{X_{2} X_{3} }} & {\overline{X_{3}^{2} }}
\end{array}\right|} \]
After substituting these into (\ref{appC_6}) and making necessary simplifications we thus obtain $\alpha _{31}^{*} =0$.

The second equality is proved in just the same way.

{\textbf{Proof} of Theorem~\ref{theorem3}
 By (\ref{appC_2}) we have $b_{1} =b_{1}^{*} $. Lemma~\ref{lemma2}  shows that $b'_{1} =b''_{1} =a_{1}^{*} \left(c_{12} ,c_{13} \right)$. In what follows we need an appropriate generalization of the theorem \citet{PV10}. It is provided below in  Appendix 4. Applying it, we get

\begin{equation} \label{appC_7}
b'_{1} =b_{1}^{*} +b_{3}^{*} \alpha _{31}^{*}
\end{equation}

This proves Theorem~\ref{theorem3}, since $\alpha _{31}^{*} =0$ by lemma~\ref{lemma3}.

\section*{Appendix 4}
\label{App_D}
\subsection*{A theorem on relationship among regression coefficients}

What follows is a statement of the theorem used in Appendix 3.

\textbf{Theorem} \citep{Panov}
 Let there be an outcome Y and a set of predictors $X_{1} ,X_{2} ,\ldots ,X_{k} $.
 Consider a multiple regression equation for the outcome on the set of predictors

\begin{equation} \label{appD_1}
y=b_{0} +\sum _{i=1}^{k}b_{i} x_{i}
\end{equation}

From the set of predictors $X_{1} ,X_{2} ,\ldots ,X_{k} $, we extract a subset $\left\{X_{i_{1} } ,X_{i_{2} }
,\ldots ,X_{i_{m} } \right\}$ and introduce regression equations for each predictor on the subset of predictors extracted

\begin{equation} \label{appD_2}
x_{i} =c_{i} +\sum _{j=1}^{m}c_{i,i_{j} } x_{i_{j} }
\end{equation}

We suppose that $c_{i} =0,c_{i,i_{j} } =\delta _{i,i_{j} } $ for $i\in \left\{i_{1} ,i_{2} ,...,i_{m} \right\}$.

Finally, let there be a multiple regression equation for outcome Y on the set of predictors $\left\{X_{i_{1} } ,X_{i_{2} },\ldots ,X_{i_{k} } \right\}$

\begin{equation} \label{appD_3}
y=a_{0} +\sum _{j=1}^{k}a_{i_{j} } x_{i_{j} }
\end{equation}

Then

\begin{equation} \label{appD_4}
a_{i_{j} } =\sum _{i=1}^{k}b_{i} c_{i,i_{j} }
\end{equation}

This theorem has been used in Appendix 3 as follows. The set of all predictors is $\left\{X_{1}^{*} ,X_{2} ,X_{3} \right\}$, the extracted set of predictors is $\left\{X_{1}^{*} ,X_{2} \right\},i_{1} =1,i_{2} =2$. Then (\ref{appD_1}) becomes the equation (see \ref{appC_2})

\[
y=b_{0}^{*} +b_{1}^{*} x_{1}^{*} +b_{2}^{*} x_{2} +b_{3}^{*} x_{3}
\]
The (\ref{appD_2}) transforms into (\ref{appC_5}), and

\begin{equation} \label{appD_5}
c_{11} =1,c_{12} =0,c_{21} =0,c_{22} =1,c_{31} =\alpha _{31}^{*} ,c_{32} =\alpha _{32}
\end{equation}
The (\ref{appD_3}) is equation (\ref{appC_3}), and

\[a_{1} =b'_{1} ,a_{2} =b'_{2} \]
Thus (\ref{appD_4}) becomes (for $a_{1} =b'_{1} $)

\[b_{1} ^{{'} } =b_{1}^{*} c_{11} +b_{2}^{*} c_{21} +b_{3}^{*} c_{31} \]
Applying (\ref{appD_5}), we obtain (\ref{appC_7})

\[b'_{1} =b_{1}^{*} +b_{3}^{*} \alpha _{31}^{*} \]

\bibliographystyle{apa}
\bibliography{Varaksin_Panov-ref}

\end{document}